\documentclass[openacc]{rsproca_new}

\jname{rspa}
\Journal{Proc R Soc A\ }

\usepackage{cleveref}
\crefname{equation}{}{}
\Crefname{equation}{Equation}{Equations}
\Crefname{figure}{Fig.}{Figs.}
\crefname{section}{\S}{\S\S}

\newif \ifmarkedup
\newcommand{\refone}[1]{\ifmarkedup \textcolor{red}{#1} \else {#1} \fi}
\newcommand{\reftwo}[1]{\ifmarkedup \textcolor{blue}{#1} \else {#1} \fi}

\newcommand{\cc}{\mathrm{c.c.}}
\newcommand{\dd}{\partial}
\newcommand{\rd}{\mathrm{d}}
\newcommand{\pd}[2]{\frac{\partial #1}{\partial #2}}
\newcommand{\td}[2]{\frac{\rd #1}{\rd #2}}

\newcommand{\eps}{\epsilon}

\newcommand{\beq}{\begin{equation}}
\newcommand{\eeq}{\end{equation}}

\newcommand{\re}{\mbox{Re}}

\DeclareMathOperator*{\argmin}{arg\,min}

\begin{document}

\title{On the independence of the slow and fast scales in multiple-scale expansions, with application to Van der Pol's equation}

\author{
Gregory Kozyreff$^{1}$
and John R. King$^{2}$}


\address{$^{1}$Optique Nonlin\'eaire Th\'eorique, Universit\'e libre de Bruxelles (U.L.B.), CP 231, Campus de la Plaine, 1050 Bruxelles, Belgium.
(orcid 0000-0002-7126-1633)
\\
$^{2}$School of Mathematical Sciences, University of Nottingham, Nottingham NG7 2RD, UK.
(orcid 0000-0002-6228-8375)
}

\subject{Applied Mathematics, Multiple-scale expansions, Exponential asymptotics, Weakly nonlinear oscillators, Hopf bifurcation}

\keywords{multiple scales, Van der Pol, divergent series, exponentially small, beyond all orders}

\corres{Gregory Kozyreff\\
\email{gregory.kozyreff@ulb.be}}

\begin{abstract}

When implementing the method of multiple scales, one is instructed to treat the slow and fast time scales as if they were independent. Despite the intuitive motivation and the effectiveness of this perturbation method, one cannot fail to notice that these two scales relate to the same unique variable, so independence can only be formal. How sensible is it, then, to split a variable asymptotically into two independent ones? We elucidate this issue with Van der Pol's equation, a simple example of a Hopf bifurcation. The discussion involves carrying the multiple-scale analysis up to arbitrarily large order and dealing with the divergent character of the resulting asymptotic series. Using the technique of optimal truncation, we re-connect the two scales. Specifically, we show that an initial translation of the fast coordinate leads to a nontrivial, exponentially small, phase shift that depends on the slow coordinate. This phase shift breaks the independence of the slow and fast scales and is found to result from the nonlinearity. Numerical simulations confirm its existence, as well as the predicted scaling. The calculation is carried out in sufficient detail to provide confidence in the generality of our result, both in its essence and in its form.

\end{abstract}

\begin{fmtext}


\end{fmtext}
\maketitle

\section{Introduction}\label{sec:intro}

\begin{figure}
\centering
\includegraphics[width=12cm]{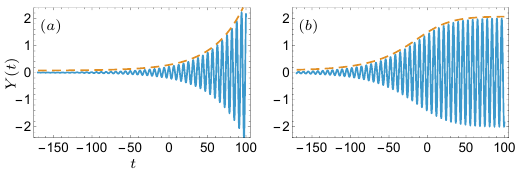}
\caption{Typical two-scale behaviours resulting from differential equations: (a) a weakly amplified linear oscillator, governed by $Y''(t)-\eps Y'(t)+Y(t)=0$ and (b) Van der Pol's equation, $Y''(t)-\eps [1-Y(t)^2] Y'(t)+Y(t)=0$ with $\eps=0.05$. The dashed lines are the slowly evolving amplitudes of oscillation, as obtained by a routine multiple-scale analysis~\cite{Kozyreff-book}.}
\label{fig:twoscales}
\end{figure}
Many autonomous differential problems produce solutions $Y(t;\eps)$ that feature  infinitely separated characteristic scales in some limit $\eps\to0$. A typical situation, depicted in \Cref{fig:twoscales}, is one in which $Y(t;\eps)$ oscillates on a fast scale $\theta(t)=\omega(\eps) t$, where $\omega(\eps)\sim1+O(\eps^2)$, with an amplitude that evolves on a slow scale $\tau=\eps t$. Here, $t$ denotes time, but it could  be any other coordinate, depending on context. The method of multiple scales (MMS) is a perturbation method that is tailored for this kind of situation, as it invites one to seek an approximation of $Y$ explicitly in the form $y(\theta,\tau;\eps)$ and provides rules to construct the dependence of $y$ on $\theta$ and on $\tau$ separately~\cite{Benderbook,Hinch1991,Kozyreff-book}. In the course of the derivation, $\theta$ and $\tau$ are treated as if they were independent variables. This is not, at first, shocking for an asymptotic-minded investigator
since $\tau$ only varies by an $O(\eps)$ amount when $t$ and $\theta$ undergo $O(1)$ changes. Hence, $\tau$ may be regarded as a ``frozen'' parameter while $\theta$ varies - at least in first approximation. The intellectual discomfort grows when one realises that the independence of $\theta$ and $\tau$ is maintained by the method up to arbitrary orders of the calculation. While the obtained approximation can be excellent from a numerical point of view, there is no \emph{a priori} justification for this fact: it is as if $Y(t;\eps)$ enjoys more symmetry than what can be expected from the autonomous character of the starting differential problem.

Specifically, in the case of time-translation invariance the existence of a solution $Y(t;\eps)$ implies that $Y(t+c_t;\eps)$ is also  solution, where $c_t$ is an arbitrary constant. From this,  distinct multiple-scale approximations can be connected via  the following group of transformations:
\begin{align}
y(\theta,\tau;\eps)&\to y(\theta+c_\theta,\tau+c_\tau;\eps),
&c_\theta&= \omega(\eps)c_t,
&c_\tau&=\eps c_t.
\label{intro:group1}
\end{align}
However, MMS implies that
\beq
y(\theta,\tau;\eps)\to y(\theta+\phi,\tau;\eps)
\label{intro:group2}
\eeq
is an additional group of allowed transformations. This indeed follows from translation invariance in $t$ and the independence of $\theta$ and $\tau$. The purpose of this paper is to invalidate and correct \cref{intro:group2} in the case of Van der Pol's equation, which is a reference model for weakly nonlinear oscillations. What we will show is that the correct transformation is
\beq
y(\theta,\tau;\eps)\to y\left(\theta+\phi+\Delta\theta\left(\tau,\phi,\eps\right),\tau+\Delta\tau\left(\tau,\phi,\eps\right);\eps\right),
\label{intro:group3}
\eeq
where $\Delta\theta$ and $\Delta\tau$ are exponentially small functions such that
\begin{align}
\Delta\theta\left(\tau,\phi,\eps\right)&\sim 4\pi\Lambda e^{- \Delta k \, \pi/\eps}\times \left\{
\begin{matrix}
0,& \tau\to -\infty,\\
 \cos ( \Delta k \, \phi), &\tau \to + \infty,
\end{matrix}
\right. 
& \Delta k&=2,
\label{intro:Deltaphi}
\end{align}
\refone{$\Lambda\approx 0.023$} is a numerical constant, and
\begin{align}
\Delta\tau\left(\tau,\phi,\eps\right)&=O\left(\eps e^{- \Delta k \, \pi/\eps}\right),
&\lim_{\tau\to\pm\infty}
\Delta\tau&= 0, 
\end{align}
Formula \cref{intro:Deltaphi} will be derived for solutions that start with an arbitrarily small amplitude and settle to a limit cycle as $t\to\infty$, as in \Cref{fig:twoscales}~\!(b). The practical consequence of \cref{intro:group3} is the following: if one modifies the small-amplitude  initial conditions in \Cref{fig:twoscales}\,(b) by a phase shift $\phi$, an additional phase shift  $\Delta\theta$ \emph{that depends on $\phi$} is found as oscillations reach the limit cycle. The fact that $\Delta\theta$ depends on $\phi$ is the signature of the fact that the slow and fast scales are not independent after all. On the other hand $\Delta\tau$ represents an exponentially small, $\phi$-dependent, translation of the front depicted by the dashed line in \Cref{fig:twoscales}\,(b).

As the calculation in this paper will show, the value of  $ \Delta k$ entering  \cref{intro:Deltaphi} is  the spacing between harmonics generated by the cubic nonlinearity in Van der Pol's equation. Starting from a leading-order solution that is a combination of $e^{i\theta}$ and $e^{-i\theta}$, the weak cubic nonlinearity enriches the solution at the following order by contributions proportional to $e^{\pm i\theta}$ and $e^{\pm3i\theta}$, and at the next order, harmonics $e^{\pm5i\theta}$ are born. Pursuing to higher orders still, all odd harmonics, with separation $\Delta k=2$ in Fourier space, are \refone{progressively} added to the solution by the cubic nonlinearity. \refone{Had we} considered a more general problem involving quadratic nonlinearities, even harmonics would be generated too. Hence, \refone{the inter-harmonic} spacing, counted in units of the fundamental frequency, would be  $ \Delta k=1$, with a significant impact on the size of $\Delta\theta$. Such a distinction between quadratic and cubic nonlinearities was already underlined in a similar calculation on weakly  nonlinear wave packets~\cite{Kozyreff2023}.  \refone{However, it is likely that the calculation of exponentially small effects would considerable increase in complexity in the presence of a quadratic nonlinearity. This is indeed what we observe if we compare the  derivations of snaking bifurcation diagrams of the Swift-Hohenberg equation in the case of quadratic \cite{Chapman2009} and cubic \cite{Dean2011} nonlinearities, which are based on the same technique as the one presented here (see also \cite{Witt2019,VillarSepulveda2025} for snaking with  Schnakenberg's and Brusselator models). Hence, in this paper, we choose to focus on a problem with purely cubic nonlinearity in order to keep technical difficulties to the minimal level.} 

\refone{The presence of the inter-harmonic spacing $\Delta k$ in  \cref{intro:Deltaphi}} implies that the connection between slow and fast scales, in autonomous differential problems, is a nonlinear effect. To see why, consider one of the simplest linear multiple-scale problem, namely $Y''(t)-\eps Y'(t)+Y(t)=0$. Its general solution is, exactly, $Y=Ae^{\tau/2}e^{i\theta}+\cc$, where $\cc$ means ``complex conjugate'', and where $\tau=\eps t$, $\theta=\sqrt{1-\eps^2/4}\,t$. By the linearity of the equation, the solution is a mere combination of exponentials. Besides, the dependences on $\tau$ and on $\theta$ are factorised. As a consequence, the symmetry \cref{intro:group2} does apply to the set of solutions. This conclusion easily extends to higher-order linear autonomous \refone{differential} systems. Indeed their solutions can be expressed by a finite combination of exponentials of the type $e^{\lambda(\eps)t}$ and the multiple-scale character of the solution usually comes from the fact that one such $\lambda(\eps)$ has a real part that becomes either infinitely larger or infinitely smaller than the imaginary part in the $\eps\to0$ limit.

Observe also that $\Delta\theta$ is proportional to $e^{-\Delta k\, \pi/\eps}$ and, hence, is smaller than any finite power of $\eps$ in the small-$\eps$ limit. Its calculation therefore requires the use of special techniques, such as asymptotic evaluation of Fourier transforms with respect to the slow scale~\cite{YangAkylas1997,Wadee1999,Nixon2017} or beyond-all-order calculations~\cite{Adams2003,Kozyreff2006,Chapman2009,Kozyreff2009,Dean2011,Dean2015,Witt2019,VillarSepulveda2025}. We will use the latter technique. Our focus on Van der Pol's equation, namely
\begin{align}
Y''(t)-\eps [1-Y(t)^2] Y'(t)+Y(t)&=0,
\label{intro:VdP}
\end{align}
is motivated by the fact that it is \refone{perhaps} the simplest settings in which to discuss exponentially small phenomena with multiple scales. If $\eps<0$, any initial condition with $|Y(0)|,|Y'(0)|\ll1$ leads to exponentially decaying oscillations with time  (we henceforth drop explicit reference to $\eps$ in the argument of $Y$). Conversely, with $\eps>0$, $Y(t)$ tends  to a limit cycle as $t\to\infty$, \textit{i.e.}, a periodic solution see \Cref{fig:twoscales}\,(b). The value $\eps=0$ is therefore a Hopf bifurcation point for this dynamical system\footnote{\Cref{intro:VdP} can be regarded, for $\eps>0$, as a rescaled form of $u''(t)-\eps u'(t)+u(t)+u^2(t)u'(t)=0$, yielding the approximate limit cycle $u\sim2\sqrt{\eps}\cos[\omega(\eps)t+\phi]$ with the familiar square root law for the amplitude.}. \refone{Note that, while $\eps$ does not multiply the highest derivative, the limit $\eps\to0$ is singular. In that limit, the real part of the eigenvalues of the equation linearised about $Y=0$ tend to zero, as in a typical Hopf bifurcation scenario.} The limit cycle can be expressed as
\begin{align}
Y_\infty(t)&\sim \sum_{j\geq0}\eps^{j}\sum_{k\text{ odd}, |k|\leq 2j+1} a_{j,k} e^{i k \theta},
&\theta=\omega(\eps)t.
\label{intro:limitcycle}
\end{align}
The coefficients $a_{jk}$ appearing above as well the terms in the expansion
\beq
\omega(\eps)\sim 1-\frac{\eps^2}{16}+\frac{17\eps^4}{3072}+\frac{35\eps^6}{884\,736}-\frac{678\,899\eps^8}{4\,096\,079\,360}+\dotsb= 1+\sum_{n\geq1}\eps^{2n}\omega_{2n}
\label{intro:omega}
\eeq
have been computed up to very high orders~\cite{Deprit1979,Andersen1982,Buonomo1998,Amore2018}. The precise numerical knowledge of $Y_\infty(t)$ will be useful to monitor the phase of oscillations in the long-time limit of solutions such as depicted in \Cref{fig:twoscales}\,(b).

The rest of the paper is organised as follows. In \cref{sec:leading}, we present the classical two-scale analysis of \cref{intro:VdP} that describes the slow evolution of small-amplitude initial conditions towards the limit cycle, assuming that $0<\eps\ll1$. In \cref{sec:late}, we derive approximate expression for the  $O(\eps^j)$ terms of the solution in the limit $j\gg1$ (the calculation is taken from  \cite{Kozyreff-book} and adapted to the present discussion). The key observation, in the spirit of Dingle's theory~\cite{Dingle1973} and of~\cite{OldeDaalhuis1995}, is that the asymptotic series generated by MMS ultimately diverges. Truncating that series near its smallest terms leaves an exponentially small remainder that connects the slow and fast scales and yields the sought-after phase shift $\Delta\theta$. Finally, in \cref{sec:num}, we check \cref{intro:Deltaphi} against numerical simulations and confirm the predicted oscillations and their exponentially small factor.

\section{Multiple-scale analysis: leading orders}\label{sec:leading}
Casting $Y(t)$ as $y(\theta,\tau)$ where $\theta$ is defined in \cref{intro:limitcycle,intro:omega}, MMS interprets $Y'(t)$ as $\omega(\eps)\dd_\theta y+\eps\dd_\tau y$, where $\dd_\theta$ and $\dd_\tau$ denote partial differentiation with respect to $\theta$ and $\tau$, respectively. The multiple-scale version of \cref{intro:VdP} is 
\beq
\omega^2(\eps)\dd^2_\theta y+\eps\omega(\eps)\dd_\theta\left(2\dd_\tau y-y\right)+\eps^2 \dd_\tau\left(\dd_\tau y- y\right)
+ \eps\left[\omega(\eps)\dd_\theta+\eps\dd_\tau\right]\left( \frac{1}{3}y^3\right)
+y =0.
\label{eq:VdP:MS}
\eeq
We start by assuming the expansion
\begin{align}
y&\sim\sum_{j=0}^\infty\eps^jy_j(\theta,\tau),
&y_j&=\sum_{k=-2j-1}^{2j+1}A_{j,k}(\tau)e^{ik(\theta+\phi)},
& (k \text{ odd})
\label{eq:VdP:series}
\end{align}
where $\phi\in\mathbb R$ is an arbitrary constant. \reftwo{In the second sum above, we have anticipated the fact that only odd harmonics are generated by the nonlinearity and that one pair of extra harmonics appear in the formal solution each time the order $j$ increases, as discussed in the introduction. We then have}
 $y^3=\sum_{j,k}\eps^jC_{j,k}(\tau) e^{ik(\theta+\phi)}$, where
\beq
C_{j,k}=\sum_{j_1,j_2,j_3}\sum_{k_1,k_2,k_3}A_{j_1,k_1}(\tau)A_{j_2,k_2}(\tau)A_{j_3,k_3}(\tau),
\label{eq:vdP:Cjk}
\eeq
and where the sums are constrained by  $j_1+j_2+j_3=j\ge0$ and  $k_1+k_2+k_3=k$. Decomposing \cref{eq:VdP:MS} into separate orders $j$ and harmonics $k$, we find
\begin{multline}
-k^2 \left(A_{j,k}+2\omega_2A_{j-2,k}+\dotsb\right)\\
+ik \left[2A'_{j-1,k}-A_{j-1,k}+\omega_2\left(2A'_{j-3,k}-A_{j-3,k}\right)+\dotsb\right]
+A_{j-2,k}''-A_{j-2,k}'
\\
+\frac13 ik\left(C_{j-1,k}+\omega_2C_{j-3,k}+\dotsb\right)
+\frac13  C'_{j-2,k}
+A_{j,k}=0.
\label{eq:vdP:Ajk}
\end{multline}
Above,  $\omega_2=-1/16$ [see \cref{intro:omega}]. For $j=0$, the recurrence is
\beq
\left(1-k^2\right) A_{0,k}(\tau)=0,
\label{resonance}
\eeq
so that only $A_{0,\pm1}(\tau)$ do not vanish. Next, the equations for $j=1$ are
\beq
\left(1-k^2\right) A_{1,k}
+ik\left(2A'_{0,k}-A_{0,k}\right) +\frac13 ikC_{0,k}=0.
\eeq
In the particular case $k=\pm1$, we have $C_{0,\pm1}=3A_{0,\pm1}^2A_{0,\mp1}$ and
\begin{align}
2A'_{0,1}&=A_{0,1}-A_{0,1}^2A_{0,-1},
&
2A'_{0,-1}&=A_{0,-1}-A_{0,1}A_{0,-1}^2,
\end{align}
with the particular solution
\beq
A_{0,\pm1}=\frac{e^{\tau/2}}{\sqrt{1+e^{\tau}}}.
\label{eq:VdP:leading}
\eeq
Hence, the leading-order approximation of $Y(t)$ is $2\cos(\theta+\phi)/\sqrt{1+e^{-\tau}}$, the envelope of which is the dashed line in \Cref{fig:twoscales}\,(b). Here, the constant phase $\phi$ of the fast oscillations appears to be free and unrelated to the slow scale. We will show in what follows that it must be corrected by the exponentially small function $\Delta\theta$ discussed in \cref{sec:intro}. 

\refone{Note for future reference that $A_{0,\pm1}$ contains an infinite set of singularities in the complex plane, given by $\tau_n=(2n+1)i\pi$.}

\section{Multiple-scale analysis: late orders}\label{sec:late}
The infinite set of equations given by \cref{eq:vdP:Ajk} can in principle be solved sequentially, one value of the calculation order $j$ at a time and $k$ running from $-2j-1$ to $2j+1$ through odd values. The complexity of the calculation, however, rapidly becomes unmanageable with increasing $j$, if one wishes to retain all the analytical details. Fortunately, at large orders $j$, the following factorial-over-power ansatz  makes $A_{j,k}$  solve \cref{eq:vdP:Ajk} asymptotically as $j\to\infty$:
\begin{align}
A_{j,k}(\tau)&\sim  \left(\frac{1}{i\chi(\tau)}\right)^{j+\alpha} \Gamma(j+\alpha)  f(\tau),
&j&\gg1. 
\label{VdP:FOP1}
\end{align}
Above, $\alpha$ is an $O(1)$ constant and $f(\tau)$ is so far arbitrary, provided it has  $O(1)$ derivatives. Notice that $|A_{j-1,k}|, |C_{j-1,k}|\ll |A_{j,k}|$ and that $A_{j-1,k}'\sim -i\chi'(\tau)A_{j,k}$ in the large-$j$ limit. Substituting \cref{VdP:FOP1} into \cref{eq:vdP:Ajk} thus yields, in first approximation,
\beq
A_{j,k}(\tau)\left[1-\left(k-\chi'(\tau)\right)^2+O\left(j^{-1}\right)\right] =0.
\eeq
\reftwo{This is formally equivalent to the resonance condition \cref{resonance} found at leading orders.}
Hence, for a given non resonant harmonic $k$, the above ansatz approximately solves the linear part of the recurrence at large orders, provided that \reftwo{$k-\chi'(\tau)$ coincides with a resonant harmonic:}
\begin{align}
\chi(\tau)&=\Delta k\left(\tau-\tau_n\right),
&k-\Delta k&=\pm1.
&(\Delta k&\neq 0)
\label{chi}
\end{align}
The constant of integration $\tau_n$ can be chosen as any one of the singularities of $A_{0,1}(\tau)$ in the complex plane, namely $\tau_n=(2n+1)i\pi$, according to \cref{eq:VdP:leading}. Since $k$ only takes odd values, $\Delta k\in\{\pm2,\pm4,\pm6,\dotsb\}$.
We thus expect, for large $j$, that $y$ has the general expansion
\begin{multline}
y\sim\text{(leading orders)}+\dotsb+
\\
\eps^{-\alpha}\sum_{j\gg1}\sum_n\sum_{\Delta k}\sum_{k=\Delta k\pm1} \left(\frac{ i\eps}{\Delta k\left(\tau_n-\tau\right)}\right)^{j+\alpha}\Gamma(j+\alpha) f_{j,n,\Delta k,k}(\tau)e^{ik(\theta+\phi)},
\label{FOP2}
\end{multline} 
 where the functions $f_{j,n,\Delta k,k}$ do not depend on $j$ in first approximation, which we can express as
\beq
f_{j,n,\Delta k,k}(\tau)\sim f_{n,\Delta k,k}(\tau)+j^{-1} g_{n,\Delta k,k}(\tau)+\dotsb.
\eeq
In what follows, we will focus on $\tau_0=i\pi$. Indeed, $\tau_0$ and $\tau_{-1}=\bar\tau_0$ are the closest singularities to the real axis and \cref{FOP2} implies that they should produce the largest effect on the solution. The contribution to the large-$j$ asymptotics of $A_{j,k}$ due to $\tau_{-1}$ will be deduced from that of $\tau_{0}$ by symmetry. By the same token, we will restrict our attention to the smallest absolute values of $\Delta k$, namely $\pm2$, corresponding to $k\in\{-3,-1,1,3\}$. 
\begin{figure}
\centering
\includegraphics[width=10cm]{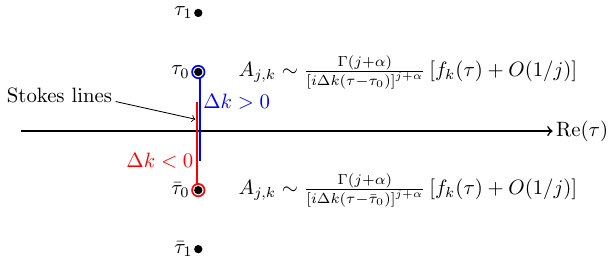}
\caption{Families of late terms and their Stokes lines in the ``slow'' complex plane. Top: terms associated to $\tau_n=\tau_0=i\pi$ in \cref{FOP2}. For such terms, $A_{j,k}$ and $A_{j+1,k}$ are of the same phase on a line that intersect the real axis (depicted in blue) provided that $\Delta k>0$. Below: terms associated to $\tau_n=\bar\tau_0$. For those, $A_{j,k}$ and $A_{j+1,k}$ are of the same phase on a line that intersect the real axis (depicted in red) provided that $\Delta k<0$. These segments are Stokes lines where exponentially small contributions to the solution are liable to be switched on \cite{Kozyreff-book,Dingle1973}.}
\label{fig:VdP:Stokes}
\end{figure}

A final consideration allows us to simplify \cref{FOP2} even further. The factorial-over-power terms appearing in this expression are known from Dingle's theory \cite{Kozyreff-book,Dingle1973} to determine Stokes lines in the complex plane where successives orders have the same phase. It is across such Stokes lines that exponentially small terms are switched on and we are thus primarily interested in those that cross the real-$\tau$ axis. Since we focus on the singularity $\tau_0$, successive terms have the same phase if $i/\Delta k(\tau_0-\tau)=1/\rho$, for some positive real number $\rho$, or equivalently, $\tau=\tau_0-i\rho/\Delta k$. Hence, the Stokes line emanating from $\tau_0$ that crosses the real axis requires $\Delta k>0$, (\Cref{fig:VdP:Stokes}). Eventually, what remains relevant to our consideration is $\Delta k=2$, corresponding to $k=1,3$. We thus have, for large order $j$, and taking only the singularity $\tau_0$ into account:
\beq
y_j\sim \sum_{k=1,3} A_{j,k} e^{ik(\theta+\phi)}
\sim\sum_{k=1,3} \left(\frac{ i}{2\left(\tau_0-\tau\right)}\right)^{j+\alpha}\Gamma(j+\alpha) \left[f_{k}(\tau) +O\left(j^{-1}\right)\right]e^{ik(\theta+\phi)},
\label{large-j}
\eeq 
%
\refone{We also need to evaluate the Fourier components $C_{j,k}$ of $y^3$ at large order $j$. While many terms make up the sum appearing in \cref{eq:vdP:Cjk}, they are dominated by a small subset of them. Given that $j\gg1$, at least one of the three indexes $j_1$, $j_2$, and $j_3$ is large. Let us assume, for, instance, that both $j_2\gg1$ and $j_3\gg1$. Then, by the factorial ansatz, $A_{j_1,k_1}A_{j_2,k_2}A_{j_3,k_3}\ll A_{j_1,k_1}A_{0,k_2}A_{j_2+j_3,k_3}$, since $\Gamma(j_2+\alpha)\Gamma(j_3+\alpha)\ll\Gamma(j_2+j_3+\alpha)$ in that limit. By the same token, wether $j_1$ is large or not, $A_{j_1,k_1}A_{0,k_2}A_{j_2+j_3,k_3}\ll A_{0,k_1}A_{0,k_2}A_{j_1+j_2+j_3,k_3}$. Recalling that the sum in  \cref{eq:vdP:Cjk} is constrained by $j_1+j_2+j_3=j$, we thus have $A_{j_1,k_1}A_{j_2,k_2}A_{j_3,k_3}\ll A_{0,k_1}A_{0,k_2}A_{j,k_3}$. We may thus asymptotically restrict the sum in  \cref{eq:vdP:Cjk} to terms where any of $j_1$, $j_2$, and $j_3$ is equal to $j$. Hence,
\begin{align}
C_{j,k}&\sim3\sum_{k_1,k_2,k_3} A_{0,k_1}A_{0,k_2}A_{j,k_3},
& (k_1+k_2+k_3=k)
\end{align}
where the factor $3$ comes from the fact that the index $j$ could arise in the first, second, or third factor of each product. Now, from the leading-order analysis (\cref{sec:leading}), $k_1$ and $k_2$ only take the values $\pm1$. The last expression thus simplifies as
\beq
C_{j,k}\sim 6 A_{0,-1}A_{0,1}A_{j,k} + 3 A_{0,-1}A_{0,-1}A_{j,k+2}+ 3 A_{0,1}A_{0,1}A_{j,k-2},
\eeq
Eventually, using \cref{large-j}, we obtain
}
\beq
C_{j,k}\sim \left(\frac{ i}{2\left(\tau_0-\tau\right)}\right)^{j+\alpha}\Gamma(j+\alpha)  \left[6A_{0,1}A_{0,-1}f_{k} +3A_{0,-1}^2f_{k+2} 
+3A_{0,1}^2f_{k-2} 
+O\left(j^{-1}\right)\right].
\eeq
Of the terms $f_k, f_{k\pm2}$ appearing above, only those equating $f_{1}$ or $f_{3}$ do not vanish. Substituting the asymptotic formulas for $A_{j,k}$ and $C_{j,k}$ we have just derived  in \cref{eq:vdP:Ajk}, all $O\left[\Gamma(j+\alpha)\right]$ contributions cancel out, leaving
\begin{align}
2f_{3}'&=f_{3}-2A_{0,1}A_{0,-1}f_{3}-A_{0,1}^2f_{1},
&
2f_{1}'&=f_{1}-2A_{0,1}A_{0,-1}f_{1}-A_{0,-1}^2f_{3},
\label{eq:f1andf3}
\end{align}
at $O\left[\Gamma(j+\alpha-1)\right]$ \refone{when we particularise  \cref{eq:vdP:Ajk} to $k=1,3$. For the other harmonics, the expressions obtained are ``only''  $O\left[\Gamma(j+\alpha-2)\right]$.} Two solutions of \cref{eq:f1andf3} are
\begin{align}
f_3(\tau)&=i\lambda_\theta A_{0,1}(\tau), 
&f_1(\tau)&=-i\lambda_\theta A_{0,-1}(\tau),
\label{beyond:late:s}
\end{align}
and
\begin{align}
f_3(\tau)&=\lambda_\tau A_{0,1}'(\tau), 
&f_1(\tau)&=\lambda_\tau A_{0,-1}'(\tau),
\label{beyond:late:tau}
\end{align}
where the coefficients $\lambda_\theta$ and $\lambda_\tau$ are arbitrary constants. They are so labelled because they are ultimately related to formal translation invariance with respect to $\theta$ and $\tau$. Indeed, if $y(\theta,\tau)$ solves the multiple-scale version of Van der Pol's equation, then so does 
\beq
y(\theta+\rd\theta,\tau)\sim y(\theta,\tau)+\rd\theta \,  \dd_\theta y(\theta,\tau)\sim  y(\theta,\tau) +\rd\theta \left[iA_{01}(\tau)e^{i (\theta+\phi) }+\cc\right]
\label{beyond:VdP:shift}
\eeq
and
\beq
y(\theta,\tau+\rd\tau)\sim y(\theta,\tau)+\rd\tau \, \dd_\tau y(\theta,\tau)\sim  y(\theta,\tau) +\rd\tau \left[A_{01}'(\tau)e^{i (\theta+\phi)}+\cc\right]. 
\eeq
We are about to find that the two solutions \cref{beyond:late:s,beyond:late:tau} switch on an exponentially small translation in $\theta$ and in $\tau$, respectively, across the Stokes line emitted by $\tau_0$. 
\refone{ Combining the two solutions that we have just identified, we find that 
\beq
A_{j,1}\sim 
-i \lambda_\theta \left(\frac{ i/2}{\tau_0-\tau}\right)^{j+\alpha_\theta}\Gamma(j+\alpha_\theta) A_{0,-1}(\tau)
+\lambda_\tau \left(\frac{ i/2}{\tau_0-\tau}\right)^{j+\alpha_\tau}\Gamma(j+\alpha_\tau) A_{0,-1}'(\tau)
\label{Aj1large}
\eeq
and
\beq
A_{j,3}\sim 
i \lambda_\theta \left(\frac{ i/2}{\tau_0-\tau}\right)^{j+\alpha_\theta}\Gamma(j+\alpha_\theta) A_{0,1}(\tau)
+\lambda_\tau \left(\frac{ i/2}{\tau_0-\tau}\right)^{j+\alpha_\tau}\Gamma(j+\alpha_\tau) A_{0,1}'(\tau)
\label{Aj3large}
\eeq}
At this stage, neither the constants $\lambda_{\theta}$, $\lambda_{\tau}$, nor the offset $\alpha$ in the factorial ansatz are determined yet. They will be when we study the expansion in the vicinity of $\tau_0$ and we will find different values, $\alpha_\theta$ and $\alpha_\tau$, of the offsets in association with $\lambda_\theta$ and $\lambda_\tau$. 

Having determined the general form of $y_j$ we note that $\eps^{j+1}y_{j+1}$ is a factor $\eps$ smaller that $\eps^{j}y_{j}$ but also approximately a factor $i(j+\alpha)/(2(\tau_0-\tau))$ larger. Assuming that the latter factor is real, the smallest term is at order $J$, such that
\begin{equation}
\frac{\eps^{J+1} y_{J+1}}{\eps^{J}y_{J}}\sim\frac{i \eps (J+\alpha)}{2\left(\tau_0-\tau\right)}\sim1.
\end{equation}
Let us examine what happens when we optimally truncate the series for $Y\sim y$ near its smallest term:
\begin{align}
Y&\sim \sum_{j=0}^J\eps^jy_j+R,
&Y^3&\sim \sum_{j=0}^J\eps^j\left(y^3\right)_j+\left(3y_0^2+\dotsb\right)R+O(R^2),
\end{align}
where $\left(y^3\right)_j$ is a short hand for $\sum_{k}C_{j,k}(\tau) e^{ik(\theta+\phi)}$. In $Y$, the two scales $\theta$ and $\tau$ are no longer treated as independent.
If we substitute this expansion in \cref{intro:VdP} and consider that all terms up to $\eps^J$ inclusive duly cancel out, then we are left with
\begin{multline}
\mathcal L R\sim-\eps\left(2\pd{^2}{\tau\dd \theta}-\pd{}{\theta}\right)\eps^Jy_J-\eps^2\left(\pd{^2}{\tau^2}-\pd{}{\tau}\right)\left(\eps^{J-1}y_{J-1}+\eps^Jy_J\right)\\
 - \frac{1}{3}\eps\left[\left(1+\eps^2\omega_2+\dotsb\right) \pd{}{\theta}+\eps\pd{}{\tau}\right]\eps^{J}\left(y^3\right)_J,
\label{eq:VdP:R1}
\end{multline}
where $\mathcal L R=\td{^2R}{t^2}-\eps\td{R}{t}+R+\frac\eps3\td{}{t}\left[\left(3y_0^2+\dotsb\right)R\right]$ is the linearised part of Van der Pol's equation. The terms in the right-hand side of \cref{eq:VdP:R1} tend to cancel each other out except in the vicinity of the Stokes lines, where $\eps^{J-1}y_{J-1}$ and $\eps^{J}y_{J}$ have the same phase. Further,   $|\eps^{J-1}y_{J-1}|$ and $|\eps^{J}y_{J}|$ only differ in relative value by an $O(J^{-1})$ quantity. Hence, in the region of interest, $\eps^{J-1}y_{J-1}\sim \eps^J y_J$. Next,  $\eps\pd{y_J}{\tau}\sim \eps (J+\alpha) y_J/(\tau_0-\tau)\sim-2i y_J$. Taking all this into account, \cref{eq:VdP:R1} simplifies, in the vicinity of the Stokes line, as
\begin{equation}
\mathcal L R\sim 4 \left(i\pd{}{ \theta}+2 \right) \eps^Jy_J + O\left(\eps^{J+1}y_J\right)
\sim 4\eps^J \left[A_{J,1} e^{i(\theta+\phi)}-A_{J,3}e^{3i(\theta+\phi)}\right].
\end{equation}
Let us now introduce the local variable $s$ as 
\beq
\tau=\tau_0-i \rho+\sqrt\eps s.
\eeq
Assuming $s=O(1)$, we may then write
\beq
\tau_0-\tau=i \rho-\sqrt\eps s
\sim i \rho e^{i\sqrt\eps s/\rho+\eps s^2/2\rho^2}
= i \rho e^{-1-i(\tau_0-\tau)/\rho+\eps s^2/2\rho^2}.
\eeq
Next, combining Stirling's formula, $\Gamma(x)\sim(2\pi/x)^{1/2}(x/e)^x$, and the condition of optimal truncation, $J+\alpha\sim 2 \rho/\eps$,  we find
\beq
\left(\frac{ i\eps}{2\left(\tau_0-\tau\right)}\right)^{J+\alpha}\Gamma(J+\alpha) 
\sim \sqrt{\frac{\pi \eps}{\rho}}e^{2i\tau_0/\eps}e^{-s^2/\rho}e^{-2i\theta}
\eeq
 in \cref{large-j}. Hence, the equation for the remainder is
\begin{equation}
\mathcal L R=4\eps^{-\alpha} \sqrt{\frac{\pi \eps}{\rho}}e^{2i\tau_0/\eps}e^{-s^2/\rho}e^{2i\phi}\left[f_1(\tau)e^{-i(\theta+\phi)}-f_3(\tau)e^{i(\theta+\phi)}\right].
\end{equation}
Note the phase factor $e^{2i\phi}$ \refone{-- we will shortly see, in \cref{mainresultminus1}, that it is responsible for the dependence of $\Delta\theta$ on $\phi$.} 
The right-hand side of the above equation for $R$ contains time through $\theta$, $\tau$ and $s$. We may thus approximately solve it by assuming the multiple-scale dependence $R=R(\theta,s,\tau)$. One readily finds
\beq
R\sim R_+=2i\pi \eps^{-\alpha} e^{2i\tau_0/\eps}e^{2i\phi} \left[f_3(\tau)e^{i(\theta+\phi)}+f_1(\tau)e^{-i(\theta+\phi)}\right]\mathcal S\left(s/\sqrt{\rho}\right),
\eeq
where  $\mathcal S$ is the Stokes function given by
\begin{align}
\mathcal S\left(x\right)&=\frac1{\sqrt\pi}\int_{-\infty}^{x}e^{-u^2}\rd u,
& \bigg(\lim_{x\to\infty}\mathcal S(x)&=1\bigg)
\label{Stokesfunction}
\end{align}
(cf.~\cite{Berry1989}) and the label $+$ refers to the singularity $\tau_0$. Finally, a contribution $R_{-}$ to the remainder comes from the singularity at $\tau=\bar\tau_0$. It must be such as to make $R$ real on the real line. Hence, we have, simply,
\beq
R\sim 2i\pi \eps^{-\alpha} e^{2i\tau_0/\eps}e^{2i\phi} \left[f_3(\tau)e^{i(\theta+\phi)}+f_1(\tau)e^{-i(\theta+\phi)}\right]\mathcal S\left(s/\sqrt{\rho}\right)+\cc.
\label{eq:VdP:R2}
\eeq
The size of $R$ is controlled by the factor \refone{$\exp (2i\tau_0/\eps)=\exp(-2\pi/\eps)$}. The  factor $2$ appearing in the exponential is the value of $\Delta k$ that was retained in the late-terms of the series. It comes from the fact that only odd harmonics are contained in the solution. Had there been a quadratic nonlinearity in Van der Pol's equation, all harmonics of the fundamental oscillation would be present and $\Delta k=1$ would be allowed, thus substantially increasing $R$.

\subsection*{Near the singularity}
It remains to determine the \refone{two} constants $\lambda_{\theta,\tau}$ and the \refone{two} offsets $\alpha_{\theta,\tau}$. \refone{This can be done by examining the behaviour of the solution in the vicinity of the singularity $\tau_0=i\pi$. As $\tau\to i\pi$, we note that $A_{0,\pm1}\sim i/\sqrt{i\pi-\tau}$. Now the equation  \cref{eq:vdP:Ajk} for $A_{j,k}$ contains the derivative of $A_{j-1,k}$ with respect to $\tau$. Hence, by induction, all Fourier components are singular, with the strength of the singularity at $\tau=i\pi$ increasing by one unit at each order. We deduce that, for \emph{any} $j$,} 
\begin{align}
A_{j,k}&\sim  \frac{B_{j,k}}{\left(i\pi-\tau\right)^{j+1/2}}, 
&\text{as}& &\tau&\to i\pi,
&\text{with}& &B_{0,\pm1}&=i.
\label{eq:VdP:A2B}
\end{align}
\refone{Substituting this ansatz into \cref{eq:vdP:Ajk}, we derive
\begin{equation}
\left(1-k^2 \right)B_{j,k}
+2ik\nu B_{j-1,k}
+\nu(\nu-1)B_{j-2,k}
+\frac13 ik D_{j-1,k}
+\frac13 \nu D_{j-2,k}
=0.
\label{eq:vdP:Bjk}
\end{equation}
where  we have introduced the short-hand notation $\nu=j-1/2$ and where
\begin{align}
D_{j,k}&=\sum_{j_1,j_2,j_3}\sum_{k_1,k_2,k_3}B_{j_1,k_1}B_{j_2,k_2}B_{j_2,k_2},
&j_1+j_2+j_3&=j,
&k_1+k_2+k_3&=k.
\label{Djk}
\end{align}
The system of equations \cref{eq:vdP:Bjk} can be regarder as an  ``inner'' problem  in the vicinity of the singularity, to be solved with the initial values
\begin{align}
B_{0,\pm1}&=i, &B_{1,\pm1}&=\pm\frac{7}{16}, &B_{1,\pm3}&=\pm\frac{1}{8}.
\end{align}
Since it is purely algebraic  \cref{eq:vdP:Bjk} is much simpler to solve iteratively than  \cref{eq:vdP:Ajk}  and the $B_{j,k}$ can be numerically evaluated up to very large $j$ using a  Computer Algebra System such as Mathematica (see Appendix and Supplementary Material for details). We should then equate \cref{eq:VdP:A2B} with the asymptotic limit as $\tau\to i\pi$ and $j\to\infty$ of  $A_{j,k}$ given by \cref{large-j}. In particular, for  $k=1,3$, \cref{Aj1large,Aj3large} yield}
\begin{align}
A_{j,1}&\sim  \left(\frac{i/2}{i\pi-\tau}\right)^{j+\alpha_s}\frac{\Gamma(j+\alpha_s)}{\left(i\pi-\tau\right)^{1/2}}  \lambda_\theta
+\frac{i}{2}
 \left(\frac{i/2}{i\pi-\tau}\right)^{j+\alpha_\tau}\frac{\Gamma(j+\alpha_\tau)}{\left(i\pi-\tau\right)^{3/2}}
\lambda_\tau , \label{Aj1}
\\
A_{j,3}&\sim  - \left(\frac{i/2}{i\pi-\tau}\right)^{j+\alpha_s}\frac{\Gamma(j+\alpha_s)}{\left(i\pi-\tau\right)^{1/2}}  \lambda_\theta
+\frac{i}{2}
 \left(\frac{i/2}{i\pi-\tau}\right)^{j+\alpha_\tau}\frac{\Gamma(j+\alpha_\tau)}{\left(i\pi-\tau\right)^{3/2}}
\lambda_\tau. \label{Aj3}
\end{align}
Matching the strengths of the singularities \refone{above with those in \cref{eq:VdP:A2B} implies}
\begin{align}
\alpha_s&=0,
&\alpha_\tau&=-1.
\end{align}
\refone{Hence, \cref{Aj1,Aj3} become
\begin{align}
A_{j,1}&\sim
 \frac{\left(i/2\right)^j \Gamma(j)}{\left(i\pi-\tau\right)^{j+1/2}}  \left(\lambda_\theta+\frac{\lambda_\tau}{j-1}\right),
&
A_{j,3}&\sim
\frac{\left(i/2\right)^j\Gamma(j)}{\left(i\pi-\tau\right)^{j+1/2}}  \left(-\lambda_\theta+\frac{\lambda_\tau}{j-1}\right).
\end{align}
Beside, we have
\begin{align}
A_{j,1}&\sim
 \frac{B_{j,1}}{\left(i\pi-\tau\right)^{j+1/2}},
&
A_{j,3}&\sim
\frac{B_{j,3}}{\left(i\pi-\tau\right)^{j+1/2}},
\end{align}
where $B_{j,1}$ and $B_{j,3}$ are computed by solving \cref{eq:vdP:Bjk}.
}
We thus obtain
\begin{align}
\lambda_\theta
&\sim \left(\frac{i}{2}\right)^{-j}  \frac{B_{j,1}-B_{j,3}}{2\Gamma(j)},
&&\refone{\lambda_\tau
\sim \left(\frac{i}{2}\right)^{-j}  \frac{B_{j,1}+B_{j,3}}{2\Gamma(j-1)},}
&j&\gg1.
\label{eq:lambdatheta}
\end{align}
\refone{Now if we make the change of variable $b_{j,k}=(i/2)^{-j}B_{j,k}$ in \cref{eq:vdP:Bjk},} we realise that the coefficients of the recurrence for $b_{j,k}$ are all real. Hence, given that $b_{j,\pm1}$ are purely imaginary, so are all subsequently computed $b_{j,k}$. We may thus assert that the right hand sides in  \cref{eq:lambdatheta} are purely imaginary. In the Appendix, we show that
\begin{align}
\lambda_\theta &= -i\Lambda,
&\Lambda&\approx 0.023080.
\label{eq:Lambda}
\end{align}
\refone{
Furthermore, we find that
\begin{align}
\lambda_\tau &= i\Gamma,
&\Gamma&\approx 0.045989.
\end{align}
}
Neglecting the contribution to $A_{j,k}$ associated with $\lambda_\tau$, which is $O(J^{-1})$ smaller and which tends to zero as $\tau\to\infty$ anyway, we have, in \cref{eq:VdP:R2},
\begin{align}
f_3(\tau)&=\Lambda A_{0,1}(\tau), 
&
f_1(\tau)&=-\Lambda A_{0,-1}(\tau). 
\end{align}
Next, as $s/\sqrt\rho\to\infty$, $\mathcal S(s/\sqrt\rho\to\infty)\to1$ and as $\tau\to\infty$, $A_{0,\pm1}\to1$. Hence, the remainder given by \cref{eq:VdP:R2} becomes
\begin{multline}
R \sim 4\pi\Lambda \re
\left\{
e^{2i\tau_0/\eps}e^{2i\phi}\left[iA_{0,1}e^{i(\theta+\phi)}-iA_{0,1}e^{-i(\theta+\phi)}\right]
\right\}
\\
=-8\pi \Lambda e^{-2\pi/\eps}\cos(2\phi)\sin\left(\theta+\phi\right)
\\
= 4\pi \Lambda e^{-2\pi/\eps}\cos(2\phi) \dd_\theta\left[2\cos\left(\theta+\phi\right)\right]
\label{mainresultminus1}
\end{multline}
Comparing this last expression with \cref{beyond:VdP:shift}, and given that $2\cos(\theta+\phi)$ is the leading-order approximation of the solution in the long-time limit, we thus find that
\beq
\lim_{\tau\to\infty}\Delta\theta \sim   4\pi \Lambda e^{-2\pi/\eps}\cos(2\phi).
\label{mainresult}
\eeq
This last expression is the main result of this paper. According to it, the maximum variation of $\Delta \theta$, 
 between $\cos(2\phi)=1$ and $\cos(2\phi)=-1$, is 
\beq
\Delta(\eps)\equiv 8\pi\Lambda \,e^{-2\pi/\eps}\approx 0.58 \,e^{-2\pi/\eps}.
\label{scaling:Delta}
\eeq

\section{Numerical investigation}\label{sec:num}
To test our prediction, we compute a large set of numerical solutions  $Y_\text{num}(t)$  of equation \cref{intro:VdP} with varying initial phases of oscillation.  Setting the initial amplitude of oscillations to the value $(1+e^{20})^{-1/2}$, the nonlinear term in \cref{intro:VdP} is completely negligible, so that $Y(t)$ is very well approximated at early stages of its evolution by 
\[
\frac{\re[e^{(\eps+i\sqrt{1-\eps^2/4})t+i\phi}]}{\sqrt{1+e^{20}}}.
\]
Hence, an initial condition with precise phase  $\phi$ is
\begin{align}
Y_\text{num}(0)&=\frac{\cos\phi}{\sqrt{1+e^{20}}}\,,
&
Y_\text{num}'(0)&=\frac{\re[(\eps+i\sqrt{1-\eps^2/4})e^{i\phi}]}{\sqrt{1+e^{20}}}\,.
\label{num:inicon}
\end{align}
At the same time, we define $Y_\infty(t)$ as being the limit-cycle solution given by \cref{intro:limitcycle} subject to $Y_\infty'(0)=0$. Since the limit cycle attracts all trajectories in the phase plane other than $Y\equiv0$, we necessarily have $Y_\text{num}(t)\to Y_\infty(t+\Delta t(\phi,\eps))$ as $t\to\infty$. The final phase of oscillations is, therefore,
\beq
\phi_\infty(\phi,\eps)=\omega(\eps) \Delta t(\phi,\eps) ,
\eeq
where $\omega(\eps)$ is given explicitly by \cref{intro:omega}. It breaks down into
\beq
\phi_\infty(\phi,\eps)=\phi+\phi_\text{dyn}(\eps)+\Delta\theta(\phi,\eps).
\eeq
Above, $\Delta\theta(\phi,\eps)$ is the quantity of interest and $\phi_\text{dyn}(\eps)$ is a dynamic phase. The latter results from the slow evolution of  the  frequency from $\sqrt{1-\eps^2/4}$  to  $\omega(\eps)$ as the amplitude ramps up towards the limit cycle, independently of $\phi$. Numerically, $\Delta t(\phi,\eps)$ is computed as
\begin{align}
\Delta t(\phi,\eps) 
&=\argmin_{s}\int_{40/\eps}^{40/\eps+2\pi}\left[Y_\infty(t)-Y_\text{num}(t-s)\right]^2\rd t,
&
s&\in[0,2\pi/\omega(\eps)].
\end{align}

\begin{figure}
\centering
\includegraphics[width=13cm]{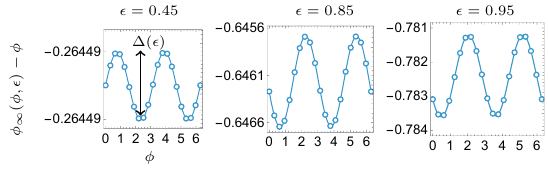}
\caption{Numerically computed phase shift of the oscillations after the establishment of the nonlinear oscillations, as a function of the initial phase, $\phi$, see \cref{num:inicon}, for three values of $\eps$. If the evolution of the complex amplitude of oscillations was independent of $\phi$, this phase shift should be constant. However it oscillates with $\phi$ in agreement with \cref{mainresult}. We define on these graphs $\Delta(\eps)$ as the maximum variation of the phase shift.}\label{figphaseshift}
\end{figure}

\begin{figure}
\centering
\includegraphics[width=7cm]{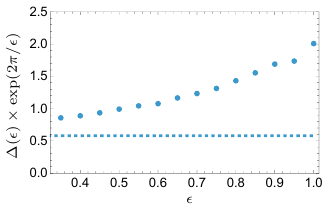}
\caption{Comparison between the analytical prediction \cref{scaling:Delta} (dashed line) and values drawn from \Cref{table:Delta} for the amplitude of oscillations of the phase shift, $\Delta(\eps)$.}\label{figscaling}
\end{figure}

In \Cref{figphaseshift}, we plot the phase shift $\phi_\infty(\phi,\eps)-\phi$ as a function of $\phi$ for three values of $\eps$ and find confirmation that the quantity oscillates as $\cos(2\phi+\text{const})$, in agreement with \cref{mainresult}. The constant in the argument of the foregoing cosine is due to a slight change of the way $t=0$ is defined in the present numerical section compared to \cref{sec:leading}. We were able to  compute the amplitude $\Delta(\eps)$ of the phase shift numerically for values of $\eps$ down to $\eps=0.35$, see \Cref{table:Delta}. A good agreement with the analytically predicted scaling is demonstrated in \Cref{figscaling}.

\begin{table}
\caption{Numerical computation of $\Delta(\eps)$ as defined in \Cref{figphaseshift}.}
\label{table:Delta}
\begin{tabular}{l c}
$\eps$& $\Delta(\eps)$\\
\hline

0.35& $1.37\times 10^{-8}$\\
0.4& $1.34\times 10^{-7}$\\
0.45& $8.07\times 10^{-7}$\\
0.5& $3.45\times 10^{-6}$\\
0.55& $1.14\times 10^{-5}$\\
0.6& $3.05\times 10^{-5}$\\
0.65& $7.39\times 10^{-5}$\\
0.7& $1.56\times 10^{-4}$\\
0.75& $3.02\times 10^{-4}$\\
0.8& $5.56\times 10^{-4}$\\
0.85& $9.58\times 10^{-4}$\\
0.9& $1.57\times 10^{-3}$\\
0.95& $2.33\times 10^{-3}$\\
1& $3.75\times 10^{-3}$\\
\end{tabular}
\end{table}

\section{Conclusion}
We have attempted to show in this paper to which extent the working assumption of scale independence is erroneous in the method of multiples scales (MMS). The interest of this perturbation method precisely lies in the fact that the two scales approximately decouple, in the sense that we can separately treat phenomena occurring on each scale, and thus significantly decrease the technical difficulty of the problem under consideration. The actual coupling between the two scales is indeed exponentially weak, as was already found in previous studies~\cite{YangAkylas1997,Wadee1999,Adams2003,Kozyreff2006,Chapman2009,Kozyreff2009,Dean2011,Dean2015,Witt2019,VillarSepulveda2025}. It is characteristic of most of these calculations  that, while they can be quite involved, the end result is often comparatively  simple. This pattern clearly repeats itself in the present setting, namely a uniformly attracting limit cycle. The end formula  \cref{mainresult} is indeed compact and easily verifiable. In the present case, the exponentially small phase shift $\Delta\theta$ may be regarded as a benign consequence of slow scale-fast scale interaction. However, in other settings, the resulting effect can be much more dramatic, notably by conditioning the existence of localised oscillatory states and ruling their snaking bifurcation diagrams~\cite{YangAkylas1997,Wadee1999,Adams2003,Kozyreff2006,Chapman2009,Kozyreff2009,Dean2011,Dean2015,Witt2019,VillarSepulveda2025} and by invalidating the classical statement that fundamental solitons of the nonlinear-Schr\"{o}dinger type travel at constant velocity~\cite{Kozyreff2023}.

The details of the calculation allows us to locate where the interaction between the slow and fast scales takes place: in the vicinity of the Stokes lines, which is also where the square of the amplitude of oscillation has its fastest growth. Well away from this Stokes lines, we know from the Stokes function \cref{Stokesfunction} that the slow and fast scales practically do not interact, even under exponentially magnifying lenses. This confirms that the underlying assumption of MMS is reliable. 

The present study also serves to make the path of exponential asymptotics with multiple-scale problems a well-beaten one. One may now indeed envisage to study other Hopf bifurcations, with slightly more complicated nonlinearities than Van der Pol's equation. As the calculation showed, the size of the exponentially small phase shifts born at the Stokes line is related to the set of harmonics generated by the nonlinearity. Including even harmonics in the picture is bound to change the factor $\exp (-2\pi/\eps)$ in \cref{mainresult} to $\exp (-\pi/\eps)$, making beyond-all-order effects numerically more visible. Technically, this results from the values $\Delta k$ allowed in the factorial-over-power ansatz \cref{FOP2}. Quadratic nonlinearities lead to the possibility that $\Delta k=\pm1$. This in turns bring the order of optimal truncation to $J\sim |\tau_0-\tau |/\eps$ instead of  $2|\tau_0-\tau| / \eps$. \refone{It would therefore be interesting to apply  the calculation of the present paper to a more general Hopf bifurcation problem involving quadratic nonlinearities.}

\ack{G.K. is a Senior Research Associate of the Fonds de la Recherche Scientifique - FNRS (Belgium.)}

\appendix
\section*{Appendix: determination of $\Lambda$}
We give here more details on the resolution of the multiple-scale recurrence in the vicinity of $\tau=i\pi$.  
\refone{By solving \cref{eq:vdP:Bjk} over a few iterations, one quickly realizes that $B_{j,-k}=(-1)^j B_{j,k}$.} For a given $j$, we thus only need to compute $B_{j,k}$ for $k=1,3,5,\dotsb,2j+1$. In so doing, we have to distinguish the harmonic $k=1$ from the others. \refone{For $3 \leq k \leq 2j+1$, once} the amplitude $B_{j-1,k}$ have been determined, we may evaluate $D_{j-1,k}$ and immediately compute
\beq
B_{j,k}=\left(k^2 -1 \right)^{-1}\left[
2ik\nu B_{j-1,k}
+\nu(\nu-1)B_{j-2,k}
+\frac13 ik D_{j-1,k}
+\frac13 \nu D_{j-2,k}\right].
\eeq
The computation of $B_{j,1}$, on the other hand, requires one to solve
\begin{equation}
2i B_{j,1}
+\nu B_{j-1,1}
+\frac{1}{3\left(\nu+1\right)} i D_{j,1}
+\frac13  D_{j-1,1}
=0,
\end{equation}
and the complication resides in the fact that $B_{j,\pm1}$ are contained in the expression of $D_{j,1}$. 

In practice, the main difficulty lies in the computation of $D_{j,k}$. Given the set of numbers $\{B_{j',k'}\}$ for $0\leq j'\leq j$, the direct implementation of \cref{Djk} turns out to be code-heavy and innefficient. Instead, we find it much more convenient to successively compute
\begin{align}
\mathcal Y_j &= \sum_{k=1}^{2j+1}\left(B_{j,k} z^{k}+B_{j,-k} z^{-k}\right),
&\left(\mathcal Y^2\right)_j &= \sum_{j'=0}^{j}\mathcal Y_{j'}\mathcal Y_{j-j'},
&\left(\mathcal Y^3\right)_j &= \sum_{j'=0}^{j}\left(\mathcal Y^2\right)_{j'}\mathcal Y_{j-j'},
\end{align}
and extract $D_{j,k}$ as the coefficient of $z^k$ in $\left(\mathcal Y^3\right)_j$. \refone{This last operation can be executed with great efficiency} with a Computer Algebra System such as Mathematica, by way of the built-in command ``Coefficient'' (a code is provided as Supplementary Material).

\begin{figure}
\centering
\includegraphics[width=13cm]{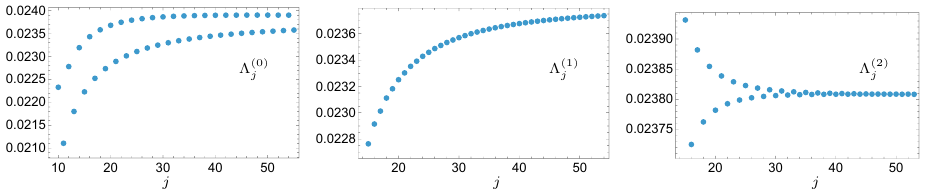}
\caption{Convergence and accelerated convergence of the series $\Lambda^{(0)}_j$ towards the value $\Lambda$.}\label{figLambda0}
\end{figure}

Once the numbers $B_{j,k}$ have been computed, we generate the sequence
\beq
\Lambda^{(0)}_j= i \left(\frac{i}{2}\right)^{-j}  \frac{B_{j,1}-B_{j,3}}{2\Gamma(j)},
\eeq
which is expected to tend to its limiting value $\Lambda$ with only an $O(j^{-2})$ error. Indeed, we found that that $O(j^{-1})$ terms, which are proportional to $\lambda_\tau$, cancel each other out  in $\Lambda^{(0)}_j$. \Cref{figLambda0} indicates that $\Lambda^{(0)}_j$ tends to its limit in an  oscillatory fashion. We first suppress this behaviour using Shanks transformation~\cite{Kozyreff-book}:
\beq
\Lambda^{(1)}_j= \frac{\Lambda^{(0)}_{j+1}\Lambda^{(0)}_{j-1}-\left(\Lambda^{(0)}_{j}\right)^2}{\Lambda^{(0)}_{j+1}+\Lambda^{(0)}_{j-1}-2\Lambda^{(0)}_j}.
\eeq
The resulting series, free from oscillations, is finally accelerated using Richardson's extrapolation, taylored to remove $O(j^{-2})$ terms~\cite{Kozyreff-book}:
\beq
\Lambda^{(2)}_j= \frac{(j+1)^2\Lambda^{(1)}_{j+1}-j^2\Lambda^{(1)}_{j}}{(j+1)^2-j^2}.
\eeq
The plot of $\Lambda^{(2)}_j$ demonstrates the efficiency of the  procedure. \refone{Eventually, we find} $\Lambda\approx 0.02380.$ 
\refone{In the same manner, the sequence
\beq
\Gamma^{(0)}_j= - i \left(\frac{i}{2}\right)^{-j}  \frac{B_{j,1}+B_{j,3}}{2\Gamma(j)}
\eeq
yields the estimation $\Gamma\approx 0.045989$ after numerically accelerating its  convergence.
}


\vskip2pc

\begin{thebibliography}{99}

\bibitem{Kozyreff-book}
Kozyreff G. 2026 {\em Practical asymptotics}.
Cambridge University Press.

\bibitem{Benderbook}
Bender CM, Orszag SA. 1999 {\em Advanced Mathematical Methods for Scientists
  and Engineers}.
Springer-Verlag.

\bibitem{Hinch1991}
Hinch EJ. 1991 {\em Perturbation methods}.
Cambridge {U}niversity {P}ress.

\bibitem{Kozyreff2023}
Kozyreff G. 2023  Speed of wave packets and the nonlinear Schr{\"{o}}dinger
  equation. {\em Phys. Rev. E} \textbf{107}.
(\href{http://dx.doi.org/10.1103/physreve.107.014219}{10.1103/physreve.107.014219})

\bibitem{Chapman2009}
Chapman SJ, Kozyreff G. 2009  Exponential asymptotics of localised patterns and
  snaking bifurcation diagrams. {\em Physica D} \textbf{238}, 319.
(\href{http://dx.doi.org/10.1016/j.physd.2008.10.005}{10.1016/j.physd.2008.10.005})

\bibitem{Dean2011}
Dean AD, Matthews PC, Cox SM, King JR. 2011  Exponential asymptotics of
  homoclinic snaking. {\em Nonlinearity} \textbf{24}, 3323.
(\href{http://dx.doi.org/10.1088/0951-7715/24/12/003}{10.1088/0951-7715/24/12/003})

\bibitem{Witt2019}
de~Witt H. 2019  Beyond all order asymptotics for homoclinic snaking in a
  {S}chnakenberg system. {\em Nonlinearity} \textbf{32}, 2667--2693.
(\href{http://dx.doi.org/10.1088/1361-6544/ab0b1d}{10.1088/1361-6544/ab0b1d})

\bibitem{VillarSepulveda2025}
Villar-Sep{\'u}lveda E. 2025  Beyond-all-order asymptotics for homoclinic
  snaking of localised patterns in reaction-transport systems.
  (\href{http://dx.doi.org/10.48550/ARXIV.2501.02698}{10.48550/ARXIV.2501.02698})

\bibitem{YangAkylas1997}
Yang TS, Akylas TR. 1997  On asymmetric gravity-capillary solitary waves. {\em
  J.~Fluid Mech.} \textbf{330}, 215.

\bibitem{Wadee1999}
Wadee MK, Bassom AP. 1999  Effects of exponentially small terms in the
  perturbation approach to localized buckling. {\em Proc.~R.~Soc.~A}
  \textbf{455}, 2351.
(\href{http://dx.doi.org/10.1098/rspa.1999.0407}{10.1098/rspa.1999.0407})

\bibitem{Nixon2017}
Nixon SD, Akylas TR, Yang J. 2017  New Aspects of Exponential Asymptotics in
  Multiple-Scale Nonlinear Wave Problems. {\em Studies in Applied Mathematics}
  \textbf{139}, 223--247.
(\href{http://dx.doi.org/10.1111/sapm.12179}{10.1111/sapm.12179})

\bibitem{Adams2003}
Adams KL, King JR, Tew RH. 2003  Beyond-all-order effects in multiple-scales
  asymptotics: travelling-wave solutions to the {K}uramoto-{S}ivashinsky
  equation. {\em J.~Eng. Math.} \textbf{45}, 197.

\bibitem{Kozyreff2006}
Kozyreff G, Chapman SJ. 2006  Asymptotics of Large Bound States of Localized
  Structures. {\em Phys. Rev. Lett.} \textbf{97}, 044502.

\bibitem{Kozyreff2009}
Kozyreff G, Tlidi M, Mussot A, Louvergneaux E, Taki M, Vladimirov AG. 2009
  Localized Beating between Dynamically Generated Frequencies. {\em Phys. Rev.
  Lett.} \textbf{102}, 043905.
(\href{http://dx.doi.org/10.1103/PhysRevLett.102.043905}{10.1103/PhysRevLett.102.043905})

\bibitem{Dean2015}
Dean AD, Matthews PC, Cox SM, King JR. 2015  Orientation-Dependent Pinning and
  Homoclinic Snaking on a Planar Lattice. {\em SIAM Journal on Applied
  Dynamical Systems} \textbf{14}, 481--521.
(\href{http://dx.doi.org/10.1137/140966897}{10.1137/140966897})

\bibitem{Deprit1979}
Deprit A, Schmidt DS. 1979  Exact coefficients of the limit cycle in van der
  Pol's equation. {\em Journal of Research of the National Bureau of standards}
  \textbf{84}, 293.

\bibitem{Andersen1982}
Andersen CM, Geer JF. 1982  Power Series Expansions for the Frequency and
  Period of the Limit Cycle of the Van Der Pol Equation. {\em SIAM Journal on
  Applied Mathematics} \textbf{42}, 678--693.
(\href{http://dx.doi.org/10.1137/0142047}{10.1137/0142047})

\bibitem{Buonomo1998}
Buonomo A. 1998  The Periodic Solution of Van Der Pol's Equation. {\em SIAM
  Journal on Applied Mathematics} \textbf{59}, 156--171.
(\href{http://dx.doi.org/10.1137/s0036139997319797}{10.1137/s0036139997319797})

\bibitem{Amore2018}
Amore P, Boyd JP, Fernandez FM. 2018  High order analysis of the limit cycle of
  the van der Pol oscillator. {\em Journal of Mathematical Physics}
  \textbf{59}.
(\href{http://dx.doi.org/10.1063/1.5016961}{10.1063/1.5016961})

\bibitem{Dingle1973}
Dingle R. 1973 {\em Asymptotic expansions: their derivation and
  interpretation}.
Academic Press Inc.

\bibitem{OldeDaalhuis1995}
Olde~Daalhuis AB, Chapman SJ, King JR, Ockendon JR, Tew RH. 1995  Stokes
  Phenomenon and Matched Asymptotic Expansions. {\em SIAM Journal on Applied
  Mathematics} \textbf{55}, 1469--1483.
(\href{http://dx.doi.org/10.1137/s0036139994261769}{10.1137/s0036139994261769})

\bibitem{Berry1989}
Berry MV. 1989  Uniform asymptotic smoothing of Stokes's discontinuities. {\em
  Proc. R. Soc. A} \textbf{422}, 7--21.

\end{thebibliography}

\end{document}